\documentclass[11pt]{article} 
\usepackage{amsmath,amssymb,amsfonts,amsthm,amscd,graphicx,psfrag,epsfig}
\usepackage{color}
\definecolor{blue}{rgb}{0,0,0.7}
\definecolor{red}{rgb}{0.75, 0, 0}
\usepackage{stmaryrd}
\usepackage[titletoc,toc]{appendix}
\usepackage{comment}

\newtheorem{theorem}{Theorem}[section]

\newtheorem{theorem-definition}[theorem]{Theorem-Definition}
\newtheorem{theorem-construction}[theorem]{Theorem-Construction}
\newtheorem{lemma-definition}[theorem]{Lemma--Definition}
\newtheorem{lemma-construction}[theorem]{Lemma--Construction}
\newtheorem{lemma}[theorem]{Lemma}
\newtheorem{proposition}[theorem]{Proposition}
\newtheorem{corollary}[theorem]{Corollary}
\newtheorem{conjecture}[theorem]{Conjecture}
\newtheorem{definition}[theorem]{Definition}

\begin{document}
\newcommand{\be}{\begin{equation}}
\newcommand{\ee}{\end{equation}}
\newcommand{\bt}{\begin{theorem}}
\newcommand{\et}{\end{theorem}}
\newcommand{\bp}{\begin{proposition}}
\newcommand{\ep}{\end{proposition}}
\newcommand{\bl}{\begin{lemma}}
\newcommand{\el}{\end{lemma}}
\newcommand{\bc}{\begin{corollary}}
\newcommand{\ec}{\end{corollary}}
\newcommand{\bd}{\begin{definition}}
\newcommand{\ed}{\end{definition}}
\newcommand{\la}{\label}
\newcommand{\Z}{{\mathbb Z}}
\newcommand{\R}{{\mathbb R}}
\newcommand{\Q}{{\mathbb Q}}
\newcommand{\C}{{\mathbb C}}
\newcommand{\lms}{\longmapsto}
\newcommand{\lra}{\longrightarrow}
\newcommand{\hra}{\hookrightarrow}
\newcommand{\ra}{\rightarrow}
\newcommand{\sgn}{\rm sgn}
\begin{titlepage}
\title{Cluster ${\cal X}$-varieties at infinity}  
\author{V.V. Fock, A. B. Goncharov}
\end{titlepage}
\date{\it To Joseph Bernstein for his 70th birthday}
\maketitle

\begin{figure}[ht]
\centerline{\epsfbox{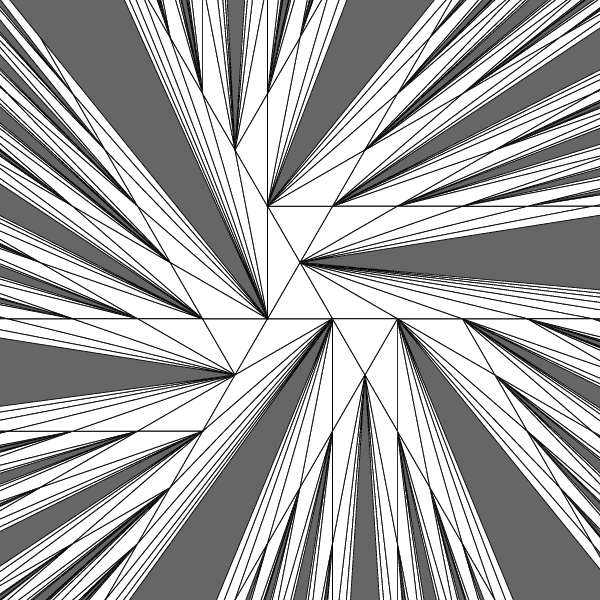}}
\caption{The tropical 
boundary hemisphere of the Teichm\"uller space  of the punctured torus.}
\label{teichmueller1}
\end{figure}

\tableofcontents

 \vskip 6mm \noindent

\begin{abstract}
A positive space is  a 
space with a positive atlas, i.e. a collection of rational coordinate systems 
with subtraction free transition functions. The set of positive real points of a positive space 
is well defined. We define  a {\it tropical compactification} of the latter. We show that 
it generalizes Thurston's compactification of a Teichm\"uller space. 

A cluster Poisson variety\footnote{originally called  cluster ${\cal X}$-variety \cite{FG2}.} 
  ${\cal X}$ is covered by a collection of 
coordinate tori $({\C}^*)^n$, which form a positive atlas of a specific kind. 
 We define a {\it special completion} $\widehat {\cal X}$ of 
${\cal X}$. It has a stratification whose strata are cluster 
Poisson varieties. The coordinate tori of ${\cal X}$ extend to coordinate affine spaces ${\Bbb A}^n$ 
in  $\widehat {\cal X}$. 

We define {\it completions of Teichm\"uller spaces} for  decorated surfaces 
${\Bbb S}$ with marked points at the boundary.
The set of positive points of the special completion of the  cluster Poisson variety 
${\cal X}_{PGL_2, {\Bbb S}}$ 
related to the Teichm\"uller theory on ${\Bbb S}$ \cite{FG1} is a part of the 
completion of the Teichm\"uller space.

\end{abstract}

\section{Introduction} 

\subsection{Geometric motivation: Thurston's boundary and 
completion  of Teichm\"uller space} 

Let $S$ be a closed oriented topological surface with $g>1$ handles. 
The Teichm\"uller space 
${\cal T}_S$ of $S$ is isomorphic to $\R^{6g-6}$. Thurston defined a compactification 
of the Teichm\"uller space by the space of {\it projective measured laminations} 
on $S$. The latter space is homeomorphic to a sphere $S^{6g-7}$. 
The action of the modular group $\Gamma_S$ of $S$ extends to the Thurston compactification. 

On the other hand, there is a completion ${\cal C}_S$ of the Teichm\"uller space, 
studied by L.Bers, H. Masur \cite{M} and others. The 
action of the modular group $\Gamma_S$ on the Teichm\"uller space extends nicely to ${\cal C}_S$. 
The quotient ${\cal C}_S/\Gamma_S$ is identified with the Deligne-Mumford 
moduli space ${\cal M}_{g}$. The components of the boundary ${\cal C}_S - {\cal T}_S$ 
are parametrised by {\it simple laminations} on $S$, that is by the isotopy classes of 
collections of pairwise non-isotopic 
non-intersecting loops on $S$. These strata are of even real codimension -- 
indeed, their quotients are complex varieties. Therefore the completion ${\cal C}_S$ 
is really different than the Thurston boundary. 

\vskip 2mm
One can extend the definition of the Teichm\"uller space and the 
Thurston compactification 
to the case when the closed surface $S$ is replaced by  a {\it decorated surface} ${\Bbb S}$, 
given by a surface $S$ 
with a finite number of holes and a finite collection of marked points on the boundary, 
considered modulo isotopy. 
Below we reserve the notation ${\Bbb S}$ for decorated surfaces, using the notation $S$ 
for the surfaces without marked points on the boudnary. 
The {\it enhanced Teichm\"uller space} ${\cal T}_{{\Bbb S}}$ for a 
decorated syrface ${\Bbb S}$ parametrises
  complex structures on the underlying surface $S$ 
plus the following additional data: marked points, and choice of   
orientations for all boundary components without marked points\footnote{When 
a hole is a puncture for  the complex structure, its boundary is empty, 
so no orientation is involved.}. 
It is diffeomorphic to an open ball, 
and the Thurston compactification is homeomorphic to a closed ball. 
The Thurston boundary  
is identified with the 
projectivisation of the space of ${\cal X}$-laminations on ${\Bbb S}$ \cite{FG3}.

We define a natural completion ${\cal C}_{\Bbb S}$ of 
the enhanced Teichm\"uller space of ${\Bbb S}$, acted by the modular group 
of ${\Bbb S}$. In the special case when
 ${\Bbb S}$ reduces to a closed surface $S$ without holes, 
we recover the special completion ${\cal C}_S$ discussed above. 
For a decorated surface ${\Bbb S}$, 
the strata of the special completion ${\cal C}_{\Bbb S}$ are parametrised by 
{\it simple  ${\cal X}$-laminations} on ${\Bbb S}$. However these strata 
are not necessarily of even codimension -- they might be of codimension one, for example. 
Since we were not able to find a reference, we elaborate this below.

\begin{figure}[ht]
\centerline{\epsfbox{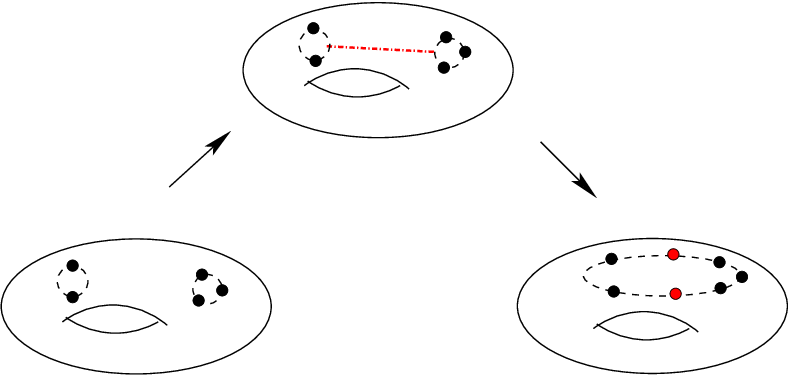}}
\caption{Starting with a torus with two holes and three marked points on the 
left, 
we cut it by a path $E$, shown on the top, getting a surface with one hole and 
five marked points.}
\label{clc0}
\end{figure}

A boundary interval is a segment of the boundary connecting two neighboring marked points. 
A simple  ${\cal X}$-lamination on ${\Bbb S}$ is a collection of pairwise non-intersecting 
loops and non-boundary, i.e. non-homotopic to a boundary interval, 
paths connecting unmarked boundary points, modulo isotopies. 
The boundary component corresponding 
to a single simple loop is a real codimension two stratum obtained 
by pinching the loop. The boundary component corresponding 
to a single non-boundary path $E$ connecting boundary points is a real codimension one stratum. 
It is identified with the enhanced Teichm\"uller space ${\cal T}_{{\Bbb S}_E}$ of the surface 
${\Bbb S}_E$ obtained by cutting the original surface along $E$, see Fig \ref{clc0}.
The marked points on ${\Bbb S}_E$ are the ones inherited from ${\Bbb S}$ 
plus the midpoints of the two components of the cutted segment $E$. 
There is a canonical fibration 
$\pi_E: {\cal T}_{{\Bbb S}} \lra {\cal T}_{{\Bbb S}_E} $ over the boundary component   
with one-dimensional fibers. More details can be found in Section 2.3.

{\bf Example}. Let $S$ be a genus $g$ oriented hyperbolic surface with $n$ holes. The classical Teichmuller space 
${\cal T}^{\rm un}_S$ 
parametrises complex structures on $S$ such that each hole is a 
puncture for the complex structure. 
Its quotient under the action of modular group of $S$ is the moduli space ${\cal M}_{g,n}$. 
The quotient of the completion of the Teichmuller space ${\cal T}^{\rm un}_S$ 
is the Deligne-Mumford moduli space $\overline {\cal M}_{g,n}$: 
a regular projective scheme over $\Z$ compactifying the moduli space ${\cal M}_{g,n}$. 
The enhanced Teichmuller space ${\cal T}_S$ is bigger: 
${\rm dim}{\cal T}_S = {\rm dim}{\cal T}^{\rm un}_S+n$.   It does not have a complex structure. 
Its completion has codimension one boundary components which do not 
intersect the completion of ${\cal T}^{\rm un}_S$. For example, if $g=n=1$, 
the classical Teichmuller space 
${\cal T}^{\rm un}_S$ is the hyperbolic plane. Its completion is obtained by adding cusps at 
the boundary of the hyperbolic plane. The completion of the enhanced Teichmuller space ${\cal T}_{{S}}$ 
is more complicated. It is discussed in Section 2.2.

\vskip 3mm

The enhanced Teichm\"uller space of an arbitrary decorated surface 
${\Bbb S}$ has an algebraic-geometric avatar: 
the moduli space ${\cal X}_{PGL_2, {\Bbb S}}$ of $PGL_2$-local systems on the 
underlying surface $S$ 
equipped with some additional structure near the boundary. 
The Teichm\"uller space is recovered as the space of its positive real points
 \cite{FG1}. 
The moduli space ${\cal X}_{PGL_2, {\Bbb S}}$ has a natural cluster Poisson variety structure. 
So one may ask whether / how one can see the Thurston compactification and components of the completion 
of the Teichm\"uller space  
in  terms of the cluster structure.  

In this paper we address these questions from a quite general algebraic 
perspective. We show that the Thurston compactification is a special case of 
a very general construction which can be applied 
to an arbitrary positive space.\footnote{The analog of the Thurston boundary 
of a positive space ${\cal X}$ was defined in Section 4 of \cite{FG1}. 
However the way it is attached to the positive real part of ${\cal X}$ -- the 
analog of the Teichm\"uller space --  was not discussed 
properly.}  Furthermore, 
the  components of the completion 
of the Teichm\"uller space corresponding to simple ${\cal X}$-laminations without loops 
can be defined in the context of arbitrary cluster Poisson varieties.

\subsection{Algebraic perspective: positive spaces and cluster Poisson varieties.} 
Many interesting varieties /
moduli spaces admit natural {\it positive atlases}. 
A positive atlas is a collection of rational coordinate systems on a 
space  
such that the transition functions between any two of the coordinate 
systems are subtraction free rational functions\footnote{These are the
 functions obtained from the coordinates by the operations of addition, multiplication and division, e.g. 
$(x+y^2)/(2x^3+xy)$. 
Precise definitions are given in Section 1.1 of \cite{FG2}. In particular, 
  the space is rational.}. 
The story originates in the classical theory of totally positive matrices. 
Further examples are discussed in  \cite{L}, \cite{FZ}, \cite{FG1}. 

\vskip 2mm
{\bf Example}. 
The  moduli space   ${\cal M}_{0,n+3}$ parametrizes  
$(n+3)$-tuples of distinct points $(x_1, ..., x_{n+3})$  
on ${\Bbb P}^1$ modulo the action of $PGL_2$. 
A cyclic order of the points gives rise to a positive atlas on 
${\cal M}_{0,n+3}$ as follows. 
Let $P_{n+3}$ be a convex polygon. 
We assign the points $x_i$ to its vertices, 
so that the order of points $x_i$ is compatible with  the 
clockwise  cyclic order of the vertices.
The coordinate systems 
are parametrized by the triangulations 
of the polygon. 
Given  a triangulation $T$, the coordinates $X^T_E$ are assigned to the diagonals $E$ of $T$. 
Namely, consider the unique rectangle formed by the sides of the polygon and diagonals of the triangulation, 
with the diagonal given by  $E$. 
Its vertices provide a cyclic configuration of four points on ${\Bbb P}^1$. 
We order them starting from a vertex of  $E$, 
getting a configuration of four points 
$(y_1, y_2, y_3, y_4)$ on ${\Bbb P}^1$, and set
\be \la{cr+}
X^T_E:= r^+(y_1, y_2, y_3, y_4) := \frac{(y_1-y_2)(y_3-y_4)}{(y_2-y_3)(y_1-y_4)}, \qquad 
r^+(\infty, -1, 0, z) = z. 
\ee
There are two possible orderings, which differ by a 
cyclic shift by two. Since the cyclic shift by one changes the cross-ratio to its inverse, 
the rational function $X^T_E$ is well defined. 
The obtained coordinate systems form a positive atlas. 
Different dihedral orders of the points lead to different positive atlases. 
For example, there are $12$ different positive atlases on ${\cal M}_{0,5}$. 
Each consists of  five coordinate 
systems. 
\vskip 2mm

A {\it positive space} is 
a space  ${X}$ together with a positive 
atlas. 
It is useful to separate the space from the atlas. 
Namely, the atlas itself determines a new scheme ${\cal X}$ over $\Z$ (perhaps non-separated). 
Let us describe the set of its complex points. 
Take a copy of the complex torus ${\rm T}_\alpha =  (\C^*)^n$ assigned 
to each coordinate system $\alpha$ of the atlas. The canonical coordinates 
$X_1, ..., X_n$ on this torus are the coordinates of the coordinate system $\alpha$. 
Gluing these tori using the transition functions between  
the coordinate systems of the atlas we get 
a scheme ${\cal X}$ ``materialising'' the atlas. 
So a positive atlas on  ${X}$ is nothing else but a birational isomorphism\footnote{We require it 
to satisfy an extra condition: the points of the torus  ${\rm T}_\alpha$ 
where the map ${\psi}$ is not defined are contained in a divisor given by 
a positive equation. So in particular the restriction of ${\psi}$ to 
the positive part $\R^n_{>0}$ 
of the torus  ${\rm T}_\alpha$ is well defined. Therefore we get 
a well defined subset $X(\R_{>0})$ of positive real points, sitting inside 
of the set $X(\R)$ positive points of $X$. Furthermore, $X(\R_{>0})$ is identified, in many different ways, with $\R_{>0}^n$.}  
\be \la{11.07.09.1}
{\psi}: {\cal X} \lra {X}.
\ee

It turns out that positive atlases appearing 
in many interesting  cases are of  very specific nature: 
the whole atlas is determined by a single coordinate system  
equipped with certain additional structure. 
The original examples were provided by {\it cluster algebras} \cite{FZI}. 
In this paper we consider such positive atlases of different type. The schemes materialising 
these atlases are called 
{\it cluster Poisson varieties}. They come with  a natural Poisson structure.

\paragraph{Cluster Poisson varieties.} 
Let us recall their definition, see Sections 1.2.1-1.2.4 of  \cite{FG2}.

\begin{definition}
A {\it quiver}  ${\bf s}$ is a datum 
$
\Bigl(\Lambda, \{e_i\}, (\ast, \ast), \{d_i\}\Bigr), 
$ where 
\begin{itemize}
\item $\Lambda$ is a lattice, i.e. a free abelian group, $\{e_i\}$ is a basis of $\Lambda$, 
and $d_i\in \Z_{>0}$ are multipliers;
\item $(\ast, \ast)$ a $\Z$-valued 
bilinear form on $\Lambda$ such that 
$\langle e_i, e_j\rangle := (e_i, e_j)d^{-1}_j$ is skew-symmetric. 
\end{itemize} 
\end{definition}
A quiver ${\bf s}$ can equivalently be described by a skewsymmetrisable matrix 
$\varepsilon_{ij}:= (e_i, e_j)$,  $i\in {\rm I}$.

\begin{definition}
A mutation of a quiver ${\bf s}$ 
in the direction of a basis vector $e_k$ is a new quiver $\widetilde {\bf s}$.
 The lattice and the form of the quiver $\widetilde {\bf s}$ 
are the same as of ${\bf s}$. 
The basis $\{\widetilde e_i\}$ of $\widetilde {\bf s}$ is defined by
\begin{equation} \label{12.12.04.2a}
\widetilde e_i := 
\left\{ \begin{array}{lll} e_i + (e_{i}, e_k)_+e_k
& \mbox{ if } &  i\not = k\\
-e_k& \mbox{ if } &  i = k.\end{array}\right.
\end{equation}
\end{definition}
The composition of mutations related to $e_k$ and 
$\widetilde e_k$ is an isomorphism of quivers.  
A composition of quiver mutations and quiver isomorphisms is 
called a {\it quiver cluster transformation}.

We assign to a quiver ${\bf s}$ 
 a split torus ${\rm T}_{\bf s}:= (\C^*)^{\rm I}$, called a {\it quiver torus},
  equipped with canonical coordinates $\{X_i\}$, $i\in {\rm I}$. 
A mutation in the direction of a basis vector $e_k$ gives rise to a birational isomorphism 
$\mu_{e_k}: {\rm T}_{\bf s} \lra {\rm T}_{\widetilde {\bf s}}$, acting on the set of coordinates 
$\{\widetilde X_i\}$ 
of ${\rm T}_{\widetilde {\bf s}}$ by 
\begin{equation} \label{f3}
\mu_{e_k}^*: \widetilde X_{i} \lms \left\{\begin{array}{lll} X_k^{-1}& \mbox{ if } & i=k, \\
    X_i(1+X_k^{-{\rm sgn}(\varepsilon_{ik})})^{-\varepsilon_{ik}} & \mbox{ if } &  i\neq k. \\
\end{array} \right.
\end{equation} 
An isomorphism of quivers gives rise to an isomorphism of the corresponding quiver tori. 
Thus a quiver cluster transformation gives rise to a birational 
isomorphism of the quiver tori. 

Given a quiver ${\bf s}$, we define the associated cluster Poisson variety as follows. Consider all quivers 
obtained from ${\bf s}$ by quiver cluster transformations. 
We assign to each of these quivers the corresponding quiver torus, and glue them using 
the corresponding birational 
isomorphisms. 

One can deduce from the Laurent Phenomenon Theorem \cite{FZL} 
 that if the form $(\ast, \ast)$ is non-degenerate,  
the regular functions on a cluster Poisson variety separate its points.

\vskip 2mm
{\bf Examples.} 

\begin{enumerate}

\item 
The cluster Poisson variety 
${\cal X}_{\Delta}$ provided by a  classical Dynkin diagrams $\Delta$, i.e. $\Delta=A_n$. 
Its matrix $\varepsilon_{ij}$ is obtained by taking the Cartan matrix 
of the Dynkin diagrams $\Delta$, and making it skewsymmetrizable 
by killing the $2$'s on the diagonal and 
changing signs  under the diagonal.  

\item The positive atlas on ${\cal M}_{0,n+3}$ assigned to a cyclic order of the points 
is a cluster Poisson variety atlas. 
Namely, let us define a quiver assigned to a triangulation $T$. The lattice 
$\Lambda$ is the free abelian group generated by the diagonals of the triangulation, with  
a basis given by the diagonals. 
The bilinear form is given by the adjacency matrix. 
Namely, two diagonals $E$ and $F$ of the triangulation are called adjacent 
if they share a vertex, and 
there are no diagonals of the triangulation between them. 
We set $\varepsilon_{EF} =0$ if $E$ and $F$ are not adjacent. 
If they are, $\varepsilon_{EF} =1$ if $E$ is before $F$ 
according to the clockwise orientation of the diagonals at the vertex $v$ shared by 
$E$ and $F$, and $\varepsilon_{EF} = - 1$ otherwise.
Mutations correspond to flips of the diagonals. One shows that the adjacency matrix 
changes under the flips according to formula (\ref{12.12.04.2a}), and that the coordinates $X^T_E$ 
mutate according to formulae (\ref{f3}). 
A zig-zag triangulation of $P_{n+3}$  provides a quiver of type $A_n$.

\item The moduli space ${\cal X}_{G,S}$ 
parametrising $G$-local systems on an oriented surface $S$ with punctures, 
equipped with a reduction of the structure group to a Borel subgroup near every puncture. 
Here $G$ is a split reductive algebraic group over $\Q$ with connected center, e.g. 
$PGL_n$ \cite{FG1}. The moduli space assigned to a decorated surface 
${\Bbb S}$ is discussed in Section 2.3. 

\end{enumerate}

So for each of these moduli spaces 
there is a cluster variety ${\cal X}$ and a 
birational isomorphism  (\ref{11.07.09.1}) describing the cluster atlas. 
However the map $\psi$ is almost never an isomorphism. 

\vskip 3mm
In Section \ref{Sec5} we investigate the following two problems:

\begin{enumerate}

\item 
{\it How one can compactify 
 the set of real positive points of an arbitrary positive space?}

\item
{\it How one can complete cluster Poisson varieties by adding some ``divisors at infinity''?}
\end{enumerate}

It turns out that these two problems are related.

\subsection{Tropical compactifications for positive spaces.} 

Let us recall a few more notions related to positive spaces. 

Let ${\Bbb F}$ be a semifield, i.e. a set with operations of addition, 
multiplication and division satisfying the usual axioms, e.g. ${\Bbb F} = \R_{>0}$. 
A more exotic example is the tropical semifield $\R^t$. 
This 
is the set ${\R}$  with the multiplication $\otimes $ 
and addition $\oplus$ given by 
$
a \otimes b := a+b, \quad a \oplus b := {\rm max}(a,b).
$ 
Similarly there are tropical semifields $\Z^t$ and $\Q^t$. 

A positive rational function provides a map of a semifield ${\Bbb F}$ 
to itself. Therefore a positive space  ${X}$ 
determines a set ${X}({\Bbb F})$ of its ${\Bbb F}$-valued points. 
As a set, it is isomorphic, in many different ways, 
to the set ${\Bbb F}^{{\rm dim}{X}}$. Namely, we assign to each coordinate system of the atlas 
the set ${\Bbb F}^{{\rm dim}{X}}$, and glue them using the coordinate 
transformations. The group  
of automorphisms of the positive space ${X}$ acts on the set ${X}({\Bbb F})$. 
If the space is a cluster Poisson variety ${\cal X}$, there is a distinguished subgroup of the latter, 
called 
the {\it cluster modular group $\Gamma_{\cal X}$ of ${\cal X}$}. 


A positive space ${\cal X}$ provides us with two spaces: 
the space of real positive points ${\cal X}(\R_{>0})$ and 
the space of real tropical points ${\cal X}(\R^t)$.

\paragraph{Example \cite{FG1}.} The set ${\cal X}_{G,S}(\R_{>0})$ of positive points 
of the moduli space ${\cal X}_{G,S}$ is the higher Teichm\"uller space assigned to the pair $(G,S)$. 
For $G=PGL_2$ we recover the enhanced Teichm\"uller space of $S$. 
The set ${\cal X}_{PGL_2,S}(\R^t)$ is the space of measured ${\cal X}$-laminations  on $S$.
\vskip 2mm

In Section \ref{Sec5} we relate the spaces ${\cal X}(\R_{>0})$ and 
${\cal X}(\R^t)$ as follows. 
We define the {\it tropical compactification} of the space 
${\cal X}(\R_{>0})$. Its boundary is called the 
{\it tropical boundary} of  ${\cal X}(\R_{>0})$. 
We show that the 
tropical boundary is canonically identified with the {\it spherical tropical space} 
${\cal S}{\cal X}(\R^t)$ -- the quotient of ${\cal X}(\R^t)$ by the natural action of the group 
$\R_{>0}^*$.

\paragraph{Example.}
We show that the tropical boundary of the positive space ${\cal X}_{PGL_2, S}$  
is the Thurston boundary of the classical Teichm\"uller space. 
Since the spherical tropical space of this positive space  is 
the space of projective measured laminations on $S$,  
our theorem 
in this case reduces to Thurston's theorem. 
\vskip 2mm

The spherical tropical space ${\cal S}{\cal X}(\R^t)$ is a sphere with
 a  piecewise linear structure preserved by the action of the automorphism
 group $\Gamma_{\cal X}$.  
It is a rather complicated object: 
Fig \ref{teichmueller1} shows a hemisphere of the  two dimensional spherical tropical space 
for  ${\cal X}_{PGL_2, S}$, where $S$ is a torus with a hole. 
The automorphism group $\Gamma_{\cal X} = PSL_2(\Z)$
 acts by  piecewise linear transformations of the hemisphere. 
So one can decompose it  into the maximal domains where the action is linear, 
as shown on Fig \ref{teichmueller1}.  We discuss this in detail  
in Section \ref{Sec6.3}. 

In Section 2.4 we introduce a convex structure on the spherical tropical space ${\cal S}{\cal X}(\R^t)$. 

\subsection{Special completions of cluster 
Poisson varieties.}

We show that the rational points the spherical tropical space  
parametrize  
a supply of  divisors 
which can appear as divisors at infinity. 
To specify a particular completion we need to choose 
a ``fan'' in the tropical space, which is reminiscent of the definition of toric varieties. 

An arbitrary cluster Poisson variety gives rise to a canonical ''fan'' in 
the tropical space ${\cal X}(\R^t)$. Namely, we assign to a coordinate system   
a cone in ${\cal X}(\R^t)$ consisting of all real tropical points which have non-negative coordinates 
in this coordinate system. 
We call it a {\it positive cone}. 

\begin{conjecture} \la{MCCC}
The interiers of the positive cones assigned to different cluster ${\cal X}$-coordinate systems are disjoint. 
\end{conjecture} 

The union of these cones  
is the {\it positive part of the tropical space}. 

The positive cone gives rise to a {\it special  completion} $\widehat {\cal X}$
of  the cluster Poisson variety. The strata of $\widehat {\cal X}$ 
 are cluster Poisson varieties.  The codimension one strata are described as follows. 
Each cluster coordinate $X_k$ gives rise to a codimension one cluster variety 
 on the boundary  $\widehat {\cal X} - {\cal X}$. 
Cluster coordinates $X_k$ and $X_k'$ produce the same boundary component if and only if 
there is a cluster transformation $\mu_{i_1} \circ \ldots \ldots \circ \mu_{i_n}$ with 
$i_1, ..., i_n \in {\rm I}-\{k\}$, which transforms $X_k$ to $X_k'$. 
This boundary component is given  by the equation $X_k=0$. The key observation is that, as it is clear from 
(\ref{f3}), a cluster Poisson transformation $\mu_j$ at $j \in {\rm I}-\{k\}$ transforms 
the equation $X_k=0$ to the one $X_k'=0$.  
This boundary cluster Poisson variety is determined by a quiver described 
 by the restriction of the matrix $\varepsilon_{ij}$ to $i,j\in {\rm I}-\{k\}$. 
The Poisson structure on  ${\cal X}$ extends to the special completion
 $\widehat {\cal X}$, 
so that the divisor at infinity is precisely the subscheme of zeros of the Poisson bivector field.

\vskip 2mm
One should not mix up tropical compactifications and 
cluster completions.  The latter are algebraic geometric objects, 
while the former 
are piecewise linear spheres. They are related however: the rational points of the 
tropical boundary parametrise a supply of divisors at infinity which can appear 
in different completions. 

\vskip 2mm

One can take points of the  special completion $\widehat {\cal X}$ 
with values in any semifield.

\paragraph{Examples.} 1. For the cluster variety  ${\cal X}_{A_n}$ 
the boundary of the special completion is a divisor combinatorially equivalent to the boundary of the 
$n$-dimensional Stasheff polytope, which we will discuss in detal below;   
the tropical boundary ${\cal S}{\cal X}_{A_n}(\R^t)$ is a simplicial sphere $S^{n-1}$ dual to the boundary of the 
$n$-dimensional Stasheff polytope; and $\widehat {\cal X}_{A_n}(\R_{>0})$ is the Stasheff polytope.
So there are two completely different ways to compactify 
the moduli space ${\cal X}_{A_n}(\R_{>0})$ of positive configurations of points on 
${\Bbb R}{\Bbb P}^1$, although topologically they both are spheres.  
We  continue discussion of this example in Section 2.3.

\vskip 2mm
2. The set of positive points of 
the  special completion of the moduli space ${\cal X}_{PGL_2, S}$ 
is a stratified space. Its codimension $k$ strata 
are the Teichm\"uller spaces for the decorated surfaces ${\Bbb S}$ 
obtained by cutting 
the original surface $S$ along $k$ non-intersecting ideal edges.

3. {\it Toric varieties.} A quiver with the trivial bilinear form is just a lattice 
with a basis. The special completion in this case is a product of ${\Bbb P}^1$'s. 
Other completions are toric varieties.

\vskip 2mm
In general 
the special completion $\widehat {\cal X}(\R_{>0})$ 
of the space ${\cal X}(\R_{>0})$ and the tropical compactification of ${\cal X}(\R_{>0})$ 
are different: the latter is a sphere, 
while the former may not be compact.

\vskip 2mm

The {positive part of the tropical boundary} in general 
fails to be a simplicial complex since some faces may be shared by infinitely many 
special simplices. 
The complement 
$$
(\mbox{the tropical boundary}) ~~ - ~~  (\mbox{the positive part of the tropical boundary}) 
$$
is an interesting object. It leads, at least in some cases, 
to generalizations of cluster coordinates. 
\vskip 2mm
For example, for the boundary of the Teichmuller space on one punctured torus, 
depicted on Fig 1,  the complement 
is a union of ``black rays'' parametrised by the ${\rm P}^1({\Bbb Q})$.
 They are the limits of the ``thin black cones'', 
which are artefacts of the resolution of the picture. The ``black rays'' are in fact 
parametrised naturally by the cusps 
on the boundary of the hyperbolic plane.  
\vskip 2mm
The following conjecture is a more precise version of Conjecture \ref{MCCC}. 
\begin{conjecture} \la{uss} The cluster complex (Definition 2.14 in \cite{FG2}) is 
naturally identified with a connected component of 
the positive part of the tropical boundary. 
\end{conjecture}

\paragraph{Moduli spaces of the Stokes data and degenerations.} 
A $G$-bundle with meromorphic connection on a Riemann surface $\Sigma$, with possibly irregular singularities, 
is determined by a {\it generalized monodromy data}, which consists of 
the monodromy and the Stokes data at irregular singularities. 
When $G=PGL_2$, 
the moduli space ${\cal M}{\cal S}(\Sigma)$ of the 
generalized monodromy data has a cluster Poisson variety structure, and is isomorphic to 
${\cal X}_{PGL_2, {\Bbb S}}$, where $S$ is the topological surface 
underlying the $\Sigma$, 
and the isotopy class of marked boundary points  on the boundary of $S$ is determined by 
the irregularity type of the connection.  
The strata of the special completion of 
generalized monodromy moduli space ${\cal M}{\cal S}(\Sigma)$
are obtained by 
 colliding the singularities of meromorphic connections on $\Sigma$. 
A similar result is valid for any reductive group $G$, and will be discussed in detail elsewhere.

\paragraph{Acknowledgments.} 
A.G. was 
supported by the  NSF grants  DMS-0653721, DMS-1059129 and DMS-1301776.
The first (repectively the final) draft of this paper was written when A.G. enjoyed the 
hospitality of IHES (Bures sur Yvette) in 2009 (repectively in 2015). 
He is grateful to IHES for the support. 

We are very grateful to the referee for correction of some of our 
errors, and  useful comments.  

\section{Tropical boundary and special completions of cluster Poisson varieties.} \la{Sec5}

\subsection{Tropical boundary of a positive space.} Let ${\cal X}$ be a positive space. 
A positive rational function on ${\cal X}$  is a rational function on ${\cal X}$ 
which in some, and hence any positive coordinate system 
can be presented as
a ratio of two Laurent polynomials with positive integral coefficients. 
Positive rational functions form a semifield denoted by  ${\Bbb Q}_+({\cal X})$. 
A positive rational function 
$F$ 
takes positive values on any positive point  
$x\in {\cal X}(\R_{>0})$. So $\log F(x)$ is defined. 
Choose a (countable) basis ${\bf B}$ in the semigroup ${\Bbb Q}_+({\cal X})$. 
The logarithms of values of functions $F\in {\bf B}$ on 
positive points  provide a map 
$$
{\cal X}(\R_{>0}) \lra \R^{\bf B}, \qquad x \lms \{\log F(x)\}_{F \in {\bf B}}.
$$
The projectivisation of this map is a map, obviously injective, 
$$
\tau: {\cal X}(\R_{>0}) \hra {\Bbb P}(\R^{\bf B}).
$$
\begin{definition} Let ${\cal X}$ be a positive space. 
The {\rm tropical compactification}  $\overline {{\cal X}^+}$ of  the space ${\cal X}(\R_{>0})$  is the closure 
of the subset 
$$
\tau({\cal X}(\R_{>0})) \subset {\Bbb P}(\R^{\bf E}).
$$ 
The {\rm tropical boundary} of the 
space ${\cal X}(\R_{>0})$ is the 
complement $\partial \overline {{\cal X}^+}:= \overline {{\cal X}^+} - {\cal X}^+$.
\end{definition}

\paragraph{Spherical tropical space.} The group $\R^*_{>0}$ acts on the set ${\cal X}(\R^t)$ of real tropical points of ${\cal X}$. 
A number  $\lambda \in \R^*_{>0}$ acts 
by multiplying all coordinates at some (and hence any) positive coordinate system by $\lambda$. 
The quotient of ${\cal X}(\R^t)$ by this action is 
the {\it spherical tropical space} ${\cal S}{\cal X}(\R^t)$. Obviously, 
$$
{\cal S}{\cal X}(\R^t) ~~
\mbox{is piecewise linear isomorphic to} ~~ S^{n-1}, \qquad n = {\rm dim}({\cal X}).
$$

\bp
There is a natural isomorphism 
$$
\eta: {\cal S}{\cal X}(\R^t) \stackrel{\sim}{\lra} \partial \overline {{\cal X}^+}.
$$
\ep

{\bf Proof}. Let us define the map $\eta$. 
Observe that
$$
\lim_{C \to \infty}\frac{\log (e^{Cx_1}+ ... + e^{Cx_n})}{C} = {\rm max}\{x_1, ..., x_n\}.
$$
Thus  the evaluation of a rational positive function 
$F$ on a tropical point $x \in {\cal X}(\R^t)$ can be interpreted as follows. 
Take a positive coordinate system. Let $(x_1, ...,  x_n)$  be the coordinates 
of $x$. Then 
\be \la{evtr1}
F^t(x) = \lim_{C \to \infty}\frac{\log F(e^{Cx_1}, ...,  e^{Cx_n})}{C}.
\ee 
It remains to notice that in the projective space 
$$
\Bigl(\frac{\log F_1(e^{Cx_1}, ... )}{C}: \frac{\log F_2(e^{Cx_1}, ...)}{C}: ... \Bigr) = 
\Bigl(\log F_1(e^{Cx_1}, ...): \log F_2(e^{Cx_1}, ...): ... \Bigr). 
$$
So we get an injective map $\eta$. Let us show that it is surjective. 
Choose a positive coordinate system. Let 
$(x_1(\varepsilon), ..., x_n(\varepsilon))$ be the coordinates of 
a family of points of ${\cal X}(\R_{>0})$ which
 has a limit under the map $\tau$ to the projective space. 
Consider its image in the subspace of the projective space ${\Bbb P}(\R^{\bf E})$ corresponding to the 
coordinate functions. Let $(x_1: ...: x_n)$ 
be the homogeneous coordinate of the limiting point, 
which we assume for a moment being defined, i.e. different then 
$(0:,...,:0)$. Then 
the limiting point is of form $\eta(x_1: ...:x_n)$. 
The proposition is proved.

\vskip 2mm
{\bf Example.} Let $S$ be a surface with punctures. 
According to \cite{FG1}, the space of projective measured laminations 
coincides with spherical tropical space for ${\cal X}_{PGL_2, S}$. 

Let ${\cal L}$ be the set of isotopy classes of simple loops on $S$. 
The Thurston compactification of the Teichm\"uller space is defined 
as the closure of the image of the map
$$
{\cal T}_S \lra {\Bbb P}(\R^{{\cal L}})
$$
of the Teichmuller space into the infinite dimensional real projective space, 
defined as follows. A point $x \in {\cal T}_S$ is assigned a point in the projective space 
whose homogeneous coordinates are the lengths of the loops in the metric on $S$ given by 
the point $x \in {\cal T}_S$.

\bp 
The Thurston boundary of the Teichm\"uller space ${\cal T}_S$ coincides with 
the tropical boundary 
of the cluster Poisson variety ${\cal X}_{PGL_2, S}$. 
\ep 

{\bf Proof}. The proposition is an easy consequence of 
the  following two facts, proved in \cite{FG1}. First, 
a loop $\alpha$ on $S$ gives rise to a positive regular function 
on the space ${\cal X}_{PGL_2, S}$ provided by the local system monodromy 
along $\alpha$. Second, let $\lambda_1, \lambda_2$ be the
 eigenvalues of the monodromy  around $\alpha$ 
for the $PGL_2(\R)$-local system assigned to a point $x \in {\cal T}_S$. 
They are defined up to a common factor. Then 
the length of the geodesic representing the isotopy class of $\alpha$ 
is $\frac{1}{2}\log |\lambda_1/\lambda_2|$.

\subsection{Tropical boundary of the enhanced Teichmuller space of a punctured torus} \la{Sec6.3}

The (enhanced) Teichmuller space for the punctured torus has dimension $3$. 
It is the space of positive points for the cluster Poisson variety structure 
described by the matrix 
\begin{equation} \label{1.3.04.1}
\varepsilon_{ij} = \left(\begin{array}{ccc}
0& 2& -2\\
-2&0&2\\ 
2&-2&0 \end{array}\right).
\end{equation}
So its tropical boundary is a sphere $S^2$ equipped with a piecewise linear action 
of the mapping class group 
$PSL_2(\Z)$. Below we describe this action,  
as well as the domains where it is linear.

\paragraph{A coordinate description.} The group $PSL_2(\Z)$ is generated by the 
standard generators $S$ and $T$, such that $S^2 =(ST)^3 = e$. The group $PSL_2(\Z)$ acts by piecewise linear transformations
of $\Z^3$: The generator $ST$ acts by $(x, y, z) \lms (y, z, x)$, 
and the generator $T$ acts by 
$$
x \lms x-2 {\rm max}(0, -z), \quad y \lms y + 2 {\rm max}(0, z), \quad z \lms -z.
$$ 
The action preserves 
$x+y+z$, and hence the half spaces $x+y+z>0$ and $x+y+z<0$. To visualize it, 
we describe its restriction to the plane $x+y+z = 1$. There is a triangle 
$$
T =  \{x,y,z ~|~ x, y, z \geq 0, x+y+z=1\}. 
$$
The $PSL_2(\Z)$-orbits of the triangle $T$ fill a dense part of the plane. 
Every side of the obtained triangles is shared by exactly two triangles. 
The complement is a countable collection of rays. Each ray has a unique 
vertex. Every vertex of the triangles is a vertex of a unique ray. 
Every vertex is shared by infinitely many triangles and just one ray.

\paragraph{A description via projective measured laminations.} 
We consider on the punctured torus $S'$ simple paths $\alpha$ without selfintersection 
of two types: 
{\it ideal arcs}, 
which start and end at the puncture, and {\it loops}. 
 
A rational lamination on $S$ is the following data:

\begin{itemize}

\item a collection of mutually non-intersecting and non-isotopic ideal arcs and 
simple loops with 
positive rational weights on the surface $S'$, up to isotopy.  

\item 
If a lamination contains an ideal arc, an orientation of a small disc containing the puncture.  
\end{itemize}

Real laminations are defined as a completion of the space of rational laminations. 

The space of projective measured laminations is identified with the sphere $S^2$ as follows. 

There is a decomposition of the sphere $S^2$ into a union of a circle and two hemispheres: 
$$
S^2 = S^1 \cup S^+ \cup S^-. 
$$

The equatorial circle $S^1$ parametrises
 projective measured laminations which are limits 
of the rational laminations supported on a single loop, as we will see shortly. 

The support of any other projective measured lamination contains an ideal arc.  
The sign of  the hemisphere $S^\pm$ 
where the point parametrising such a  projective measured lamination is located 
is determined by the orientation of a disc containing the puncture: 
it is $+$ if it is induced by the orientation of $S'$. 

The triangles parametrise projective laminations 
$$
a_1\alpha_1 + a_2\alpha_2 + a_3\alpha_3, \qquad a_1, a_2, a_3 \geq 0, 
$$
where the ideal arcs $\alpha_1, \alpha_2, \alpha_3$
form an ideal triangulation of the punctured torus $S'$,
 and the triple $(a_1, a_2, a_3)$ is defined up to multiplication by a positive number. 
The edges of the triangles 
parametrise projective laminations of the form $a_1\alpha_1 + a_2\alpha_2$, 
and the vertices correspond to the laminations 
given by the ideal arcs $\alpha$ on $S'$. 

Rays parametrise projective laminations 
$$
a\alpha + b\beta, \qquad a, b \geq 0, 
$$
where $\alpha$ is an ideal  arc, 
$\beta$ is a loop, and the pair $(a, b)$ is defined up to multiplication by a positive number. 
The   projective lamination $\beta$ corresponds to  the 
intersection of the ray with the equatorial circle. The   projective lamination $\alpha$ 
is the vertex of the ray shared by the triangles. 

The equatorial circle $S^1$ parametrises
 projective measured laminations  on the 
torus compactifying $S'$. 
Presenting a compact torus as  $\R^2/\Z^2$, 
the equatorial circle $S^1$ is identified with $\R^2/\R^*_{>0}$. A projective 
lamination is obtained by projecting a one dimensional subspace in $\R^2$ 
onto the torus. Integral laminations correspond to the lines with rational slopes, i.e.
 to the rational points $p/q$ on $S^1$. The action of the modular group 
$PSL_2(\Z)$ on the projective measured laminations is the action by the 
fractional linear transformations of the equatorial circle $S^1$. 

Let us compare the $S^+\cup S^1$ part of the tropical boundary on Fig \ref{teichmueller1} 
with the classical modular triangulation 
of the hyperbolic plane. There is a natural map from the former to the latter. It  
shrinks the rays onto the cusps, and is an isomorphism on the triangles.

\subsection{Special completions of cluster varieties} \la{Sec6.2}

\paragraph{Rational tropical points and strict valuations.} A valuation on a ring ${\Bbb L}$ is a map $v: {\Bbb L} \lra \Z$ such that
$$
v(L_1L_2) = v(L_1) + v(L_2), \qquad  v(L_1 + L_2) \leq {\rm min}(v(L_1), v(L_2)). 
$$
A {\it strict valuation} on a semiring 
${\Bbb L}_+$ is a map $v: {\Bbb L}_+ \lra \Z$ such that 
$$
v(L_1L_2) = v(L_1) + v(L_2), \qquad  v(L_1 + L_2) = {\rm min}(v(L_1), v(L_2)). 
$$
A valuation $v$ gives rise to 
the {\it valuation semiring}
$$
{\cal O}^+_v \subset {\Bbb L}_+, \qquad {\cal O}_v = 
\{F \in {\Bbb L}_+ ~|~ v(F) \geq 0\},
$$
and the {\it valuation ideal}
$$
{\cal M}^+_v \subset {\cal O}^+_v, \qquad {\cal M}_v = 
\{F \in {\Bbb L}_+ ~|~ v(F) > 0\}.
$$
So there is the quotient semiring
$$
{\cal R}^+_v := {\cal O}^+_v /{\cal M}^+_v. 
$$

The following lemma is an immediate consequence of the definitions.
\begin{lemma} \label{12.6.04.10v} Let ${\cal X}$ be a positive space. 
A point 
$x\in {\cal X}(\Z^t)$ provides a strict valuation $v_x$ of the semiring ${\Bbb L}_+({\cal X})$ 
given by the negative of the value of the tropicalization of a function $F$ at $x$:
\begin{equation} \label{12.6.04.5val}
v_x(F):= -F^t(x). 
\end{equation} 
\end{lemma} 
\vskip 2mm

{\bf Remark}. Let ${\cal X}$ and ${\cal A}^\vee$ be a Langlands dual 
pair of cluster varieties -- see Section 4 of \cite{FG2}. 
Then the intersection pairing ${\cal I}_{\cal A}(l,m): {\cal A}^\vee(\Z^t) \times {\cal X}(\Z^t) \lra \Z$ 
whose existence was conjectured in Conjecture 4.3 in {\it loc. cit.} 
is nothing else but the value of the 
valuation $v_m$ provided by an integral tropical point $m\in {\cal X}(\Z^t)$ on 
the positive regular function ${\Bbb I}_{\cal A}(l)$ on the ${\cal X}$-space corresponding by 
Conjecture 4.1 in {\it loc. cit.} to the 
integral tropical point $l\in {\cal A}^\vee(\Z^t)$:
$$
{\cal I}_{\cal A}(l,m) = v_m({\Bbb I}_{\cal A}(l)).
$$
There is, of course, a similar picture  when we interchange ${\cal X}$ and ${\cal A}$ spaces.

\paragraph{Divisors at infinity corresponding to the points of ${\cal S}{\cal X}(\Q^t)$.} 
Let us assume that the ring ${\Bbb L}({\cal X})$ 
is spanned by the semiring ${\Bbb L}_+({\cal X})$, and 
the strict valuation $v_x$ on ${\Bbb L}_+({\cal X})$ can be extended to a valuation on 
${\Bbb L}({\cal X})$, also denoted by $v_x$. 
So there are the 
valuation ring ${\cal O}_x$, the valuation ideal ${\cal M}_x$, and the quotient ring 
${\cal R}_x$ for the valuation $v_x$. 
Then, given $x\in {\cal X}(\Z^t)$, there is a scheme 
$$
{\cal D}_x:= {\rm Spec}({\cal R}_x). 
$$
We think about it as of 
a divisor at infinity for the space ${\cal X}$ corresponding to $x$. 
The integer $v_x(F)$ is the order of zero of the function $F$ at this divisor. 
The scheme ${\cal D}_x$ depends only on the image of 
$x$ in ${\cal S}{\cal X}(\Q^t)$.

If ${\cal X}$ has a cluster variety, ${\cal D}_x$ does not always 
have a cluster variety structure. Here is an example which does not require any {\it a priori} assumptions, 
 where this is the case.

\paragraph{Special completions of cluster varieties.} 
Given a quiver ${\bf i}$, any subset of its basis vectors 
 determines a quiver ${\bf j}$, called a {\it subquiver of the quiver ${\bf i}$}. Its lattice 
$\Lambda_{\bf j}$ is the sublattice of $\Lambda$ spanned by the 
specified basis vectors; the bilinear form on the lattice $\Lambda_{\bf j}$ is induced by the 
bilinear form on $\Lambda$. 
\begin{definition}
Let ${\cal X}$ be a cluster Poisson variety. 
A {\rm special cone}  is a cone
$C\subset {\cal X}(\R^t)$ such that 
there exists a quiver ${\bf i}$ and a subquiver ${\bf j} \subset {\bf i}$ 
so that $C$ is given by the following conditions on the tropical cluster coordinates $\{x_i\}$ 
for the quiver ${\bf i}$:
$$
C= \{l \in {\cal X}(\R^t)~|~ x_i(l) \geq  0, \quad x_c(l) = 0 ~\mbox{if} ~c \in 
{\bf j}\}. 
$$
\end{definition}
So $C$ is a cone of codimension $ |{\bf j}|$.

\vskip 2mm 
{\bf Remark}. A {\rm special simplex} is the projection of a special cone to the spherical tropical space. 
Special simplices may not form a  locally finite subset of the 
spherical tropical space. An example is given in Section \ref{Sec6.3}, where 
special triangles are the orbits of the triangle $T$.  

\vskip 2mm 

\bl
A  special cone $C\subset {\cal X}(\R^t)$ 
determines a cluster Poisson variety ${\cal X}_{C}$ of codimension 
$
{\rm dim} C.
$ 
\el

{\bf Proof}. Mutating basis 
vectors of the subquiver ${\bf j}$  
we do not change coordinates of a point $l \in C$. Indeed, the tropical coordinates $x_i$ change 
under 
the mutation in the direction $e_j$ as follows:
$$
x'_i = x_i \pm \varepsilon_{ij}{\rm max}(0, \mp x_j), \qquad i \not = j, \qquad x_j'=-x_j.
$$
So $x_j=0$ implies $x_i'=x_i$ for all $i$. 
Therefore  
the equivalence class of the quiver ${\bf j}$ is well defined by the
 special cone $C$.  The lemma is proved. 

\vskip 2mm
\bp \la{eexx} Assume that regular functions 
on a cluster Poisson variety ${\cal X}$ separate the points. 
Then there is a unique scheme over $\Z$, perhaps non-separated, called 
the {\rm special completion $\widehat {\cal X}$} of ${\cal X}$, such that

\begin{enumerate}

\item  $\widehat {\cal X}$ has a stratification whose strata 
are the affine closures of cluster varieties ${\cal X}_C$ parametrised by 
the special cones $C\subset {\cal X}(\R^t)$. 

\item
We have ${\cal X}_C \subset \widehat {\cal X}_{C'}$ if and only if 
$C' \subset  C$. 
\end{enumerate}
\ep

{\bf Proof}. 
A quiver ${\bf i}$ defines a cluster quiver torus ${\cal X}_{\bf i}$ which 
is canonically identified with 
${\Bbb G}_m^{{\rm I}}$ by the cluster coordinates. It has  a canonical 
embedding into the affine space
\be \la{sdf}
{\cal X}_{\bf i} \stackrel{\sim}{=} {\Bbb G}_m^{{\rm I}} \hra {\Bbb A}^{\rm I}
\ee
A subquiver 
${\bf j} \subset {\bf i}$ determines a cluster torus ${\cal X}_{\bf j}$ which we glue to 
the torus ${\cal X}_{\bf i}$ as follows: the torus ${\cal X}_{\bf j}$ is realized as a subvariety 
of ${\Bbb A}^{\rm I}$ in (\ref{sdf}) given by the conditions 
\be \la{the con}
X_j = 0 ~\mbox{ if $j \in {\bf i} - {\bf j}$, and $X_j \not = 0$ otherwise}.
\ee 
The key point is that mutations in the directions of basis vectors of ${\bf j}$ 
do not change condition (\ref{the con}). Therefore it provides a boundary stratum 
identified with the cluster variety corresponding to the quiver ${\bf j}$. 
This procedure is transitive: given a triple of quivers ${\bf k} \subset {\bf j}\subset {\bf i}$, 
the stratum for the subquiver  ${\bf k}\subset  {\bf i}$ is the same as the stratum 
for the subquiver  ${\bf k}\subset  {\bf j}$ in the stratum for ${\bf j}\subset {\bf i}$.
The proposition follows.

\vskip 2mm
{\bf Remark}. There is a version of Proposition \ref{eexx} where we do not take affine closures in 1). 
It does not require that  
regular functions 
on ${\cal X}$ separate the points.

\vskip 2mm
Taking the positive points of the boundary strata we get a completion $\widehat {\cal X}(\R_{>0})$ 
of the space 
${\cal X}(\R_{>0})$. It consists of cells, each isomorphic to $\R^{k}$, parametrised by the 
special cones, where $k$ is the codimension of the cone. 

For 
a cluster Poisson variety of finite type the completion $\widehat {\cal X}(\R_{>0})$ is homeomorphic 
to a disc. Its boundary is a polyhedron dual to the 
generalized Stasheff polyhedron of the Dynkin diagram. Otherwise the special completion 
is not  a manifold.

\vskip 2mm
Here are some examples.

\paragraph{Special completion of the moduli space of framed $PGL_2$-local systems.} 
Let ${\Bbb S}$ be a decorated surface. 
Let 
$\partial {\Bbb S}$ be its punctured boundary, i.e. 
the boundary minus  the marked points. The moduli space ${\cal X}_{PGL_2, {\Bbb S}}$ 
parametrised pairs $({\cal L}, \beta)$, called  
framed $PGL_2$-local systems, where ${\cal L}$ is a $PGL_2$-local system on ${\Bbb S}$, and 
$\beta$ is a flat section of the associated flat 
${\Bbb P}^1$-bundle  on $\partial {\Bbb S}$. 
It amounts to a collection of flat sections on connected components of the boundary. 
In the absence of  marked points we get 
the moduli space described in  the Example 3 at Section 1.2.

Choose a point inside of each  connected component  of $\partial {\Bbb S}$. 
An ideal triangulations of ${\Bbb S}$ 
is a triangulation of the surface with vertices at these points. 
\begin{figure}[ht]
\centerline{\epsfbox{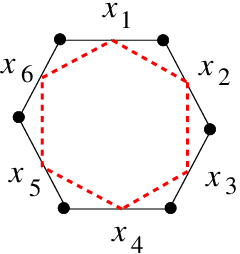}}
\caption{The convex hull of the centers of sides of the polygon $P_6$ is 
the dual polygon $\ast P_6$. A triangulation of the latter is an ideal triangulation producing a quiver.}
\label{cluste31}
\end{figure}
The moduli space ${\cal X}_{PGL_2, {\Bbb S}}$ has a cluster 
Poisson variety structure with the quivers parametrised by ideal triangulations $T$ of ${\Bbb S}$. 
The cluster coordinates for a given triangulation $T$ are assigned to the 
non-boundary edges of $T$ \cite{FG1}, just like in the case of ${\cal M}_{0, n+3}$. 

\vskip 3mm
For example, when ${\Bbb S}$ is a disc with $n+3$ marked points at the boundary, 
the moduli space ${\cal X}_{PGL_2, {\Bbb S}}$ parametrises $(n+3)$-tuples of points 
$(x_1, ..., x_{n+3})$ on ${\Bbb P}^1$ 
modulo the diagonal action of the group $PGL_2$.
To decribe its cluster structure we start with a convex $(n+3)$-gon $P_{n+3}$. 
Its vertices are the marked points. We assign the 
points $(x_1, ..., x_{n+3})$ to the sides  of the polygon $P_{n+3}$ as on Fig \ref{cluste31}.
Now take the convex hull of the centers of the sides of the polygon $P_{n+3}$. 
We get a dual convex polygon, denoted by $\ast P_{n+3}$. 
The configuration of points $(x_1, ..., x_{n+3})$ sits at its vertices. 
A triangulation of the dual  polygon $\ast P_{n+3}$ is an ideal triangulation. 

\begin{figure}[ht]
\centerline{\epsfbox{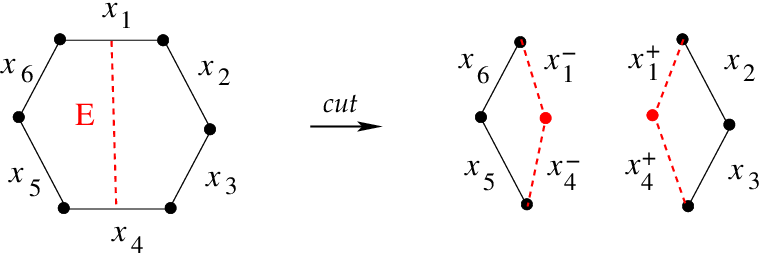}}
\caption{A boundary divisor of the cluster 
completion of ${\cal X}_{A_3}$ isomorphic to ${\cal X}_{A_1} \times {\cal X}_{A_1}$ corresponding to 
a cut of the hexagon by punctured segment, shown on the left.}
\label{cluste3}
\end{figure}

To get a boundary divisor, we cut the polygon $P_{n+3}$ along a  line segment $E$ connecting 
midponits of two of its non-consequitive sides. So $E$ is  a diagonal of the 
polygon $\ast P_{n+3}$. 
Then we add a marked point inside of  each of the two new sides of  obtained polygons. 
See Fig \ref{cluste3}. For example, cutting along an edge with the vertices 
labeled by  points $x_1, x_k$, where $k \not = n+3, 2$, we get a pair of polygons 
decorated by configurations of points 
\be \la{DEG}
(x^-_1, x_2, ..., x_{k-1}, x_k^-) ~~\mbox{and}~~(x^+_{k}, x_{k+1}, ..., x_{n+3}, x_1^+).
\ee

For a general $\Bbb S$ we choose a non-boundary 
edge $E$ of an ideal triangulation $T$ of ${\Bbb S}$. 
Let ${\Bbb S}_E$ be the decorated surface  obtained by cutting the surface along the edge $E$. 
We denote by $E_-$ and $E_+$ the pair of boundary segments on ${\Bbb S}_E$ which glue back 
to  $E$. 
The midpoints of $E_-$ and $E_+$ are new marked points. 
 The rest of the marked points are inherited 
from ${\Bbb S}$.  So the number of marked points always 
increases by two. 

There is a canonical fibration with one dimensional fibers:
\be \la{pie}
\pi_E: {\cal X}_{PGL_2, {\Bbb S}} \lra {\cal X}_{PGL_2, {\Bbb S}_E}.
\ee
Namely, given a framed local system $({\cal L}, \beta)$ on ${\Bbb S}$, 
we restrict the ${\cal L}$  to  ${\Bbb S}_E$, getting 
a local system ${\cal L}_E$. Then the framing $\beta$ 
on ${\cal L}$  induces a framing $\beta_E$ on ${\cal L}_E$. 

\begin{figure}[ht]
\centerline{\epsfbox{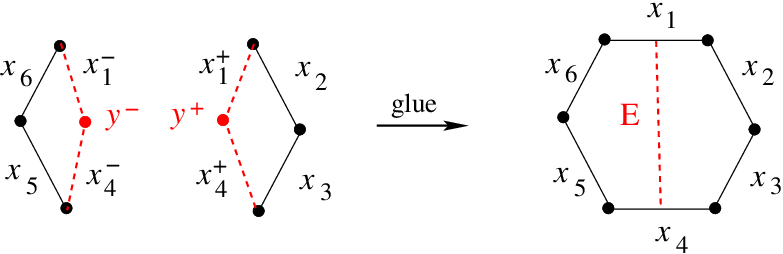}}
\caption{We add points $y^-$ and $y^+$  at the red vertices on the left, and glue the two configurations 
of points in a unique way so that the triples $(x_1^-, y^-, x_4^-)$ and $(x_1^+, y^+, x_4^+)$ are identified. 
After that we discard the $y$-points. 
Varying one of them, say $y^+$, we get the  fibers of  projection (\ref{pie}).}
\label{cluste32}
\end{figure}

The fibers of  fibration (\ref{pie}) can be parametrised, although non-canonically, as follows, 
see Fig \ref{cluste32}.
Let $({\cal L}_E, \beta_E)$ be a framed $PGL_2$-local system on 
${\Bbb S}_E$. 
Choose additional flat sections $y_\pm$  of the ${\Bbb P}^1$-bundle near the new marked points
over the boundary intervals $E_\pm$. Assume that they differ from the ones provided by the framing 
 $\beta_E$ 
on the left and on the right. Then there is a unique 
framed local system on ${\Bbb S}$ obtained from $({\cal L}_E, \beta_E)$ 
by gluing $E_-$ to $E_+$ so that the section $y_-$ is identified with $y_+$. 
Keeping $y_-$ fixed and changing $y_+$ we get the fiber. 
When $y_+$ approaches to the flat section at the end of $E_+$, we get the limiting 
framed local system on ${\Bbb S}_E$. 
He we use the orientation of the surface to tell which end.  

Special cones in the tropical space are 
parametrised by collections $\{E_1, ..., E_k\}$ of pairwise non-intersecting ideal non-boundary 
edges  on ${\Bbb S}$. 
 By the general construction, such a collection 
gives rise to a 
codimension $k$ stratum of the special completion of ${\cal X}_{PGL_2, {\Bbb S}}$. 
It is described as follows. 
Let ${\Bbb S}_{E_1, ..., E_k}$ be the surface  obtained by cutting  ${\Bbb S}$ along the edges $E_1, ..., E_k$. 
The vertices of the edges $E_1, ..., E_k$ become special points on the cutted surface. 

\bl The stratum assigned to the collection $\{E_1, ..., E_k\}$
 is naturally identified with the moduli spaces of framed $PGL_2$-local systems 
on  ${\Bbb S}_{E_1, ..., E_k}$.
\el 

{\bf Proof}. Let $T$ be an ideal triangulation  of ${\Bbb S}$. 
Denote by $Ed(T)$ the set of non-boundary edges of $T$. 
Then $Ed(T) - \{E_1, ..., E_k\}$ is the set of non-boundary edges 
of the ideal triangulation of ${\Bbb S}_{E_1, ..., E_k}$ provided by $T$. 
Let $X_{E}^T$ be the cluster coordinate for the triangulation $T$ assigned to a non-boundary 
 edge $E$ of $T$. 
Then the boundary stratum assigned to the collection $\{E_1, ..., E_k\}$ is obtained 
by letting the coordinates $X_{E}^T$, where $E \in \{E_1, ..., E_k\}$, to be zero.

\paragraph{Cluster Poisson variety ${\cal X}_{A_n}$, the Stasheff polytope and the Stasheff divisor.} The $n$-dimensional {\it Stasheff polytope} 
is a convex polytope whose vertices are parametrised 
by complete triangulations of a convex polygon $P_{n+3}$. The $k$-dimensional   
faces are parametrised by triangulations 
missing $k$ diagonals. Here is a geometric realization of the Stascheff polytope. 

A dihedral order of  points 
$(x_1, ..., x_{n+3})$ on the real projective line 
provides a connected component ${\cal M}^0_{0,n+3}(\R)$ 
of the moduli space ${\cal M}_{0,n+3}(\R)$, parametrising configurations of points on 
${\Bbb P}^1(\R)$ whose  cyclic order is compatible 
with an orientation of ${\Bbb P}^1(\R)$. 
It is nothing else but the Teichmuller space 
assigned to a disc with $n+3$ marked points on the boundary, discussed in Section 1.1. 
Its closure is 
identified with the Stasheff polytope. 

The Zariski closure of the boundary of ${\cal M}^0_{0,n+3}(\R)$ 
is the {\it Stasheff divisor} $B_n$ 
in $\overline {\cal M}_{0,n+3}$ \cite{GM}. 
Its components intersect the Stasheff polytope  
by real 
codimension one faces. 

To clarify the structure of the Stasheff divisor, recall
 some basics about $\overline {\cal M}_{0,n+3}$. 
The complement $\overline {\cal M}_{0,n+3} - {\cal M}_{0,n+3}$ 
is a normal crossing divisor whose components  
are called the {\it boundary divisors}. They 
are parametrised by decompositions
\be \la{dec}
\{x_1, ..., x_{n+3}\} = I \cup J, \qquad |{I}|, |J| >1.
\ee
Recall that we assign the points $x_1, ..., x_{n+3}$ to the vertices $p_1, ..., p_{n+3}$ 
of the convex polygon $P_{n+3}$. 

The boundary divisor 
$D_{{I},J}$, see (\ref{dec}),   belongs to the Stasheff divisor if and only if 
the convex hull of the vertices from the set $I$ does not intersect the one for the set $J$. 
So the faces of the Stasheff divisor are products of the moduli 
spaces $\overline {\cal M}_{0,m}$. 
Their structure match the combinatorial structure of the Stasheff polytope. 

\paragraph{Examples.} 1. The cross-ratio 
(\ref{cr+}) provides isomorphisms 
$$
{\cal M}_{0,4}  \stackrel{\sim}{\lra} {\Bbb P}^1 - \{\infty, -1, 0\}, \qquad 
B_1 \stackrel{\sim}{\lra} \{0, \infty\}, \qquad
\overline {\cal M}_{0,4} \stackrel{\sim}{\lra} {\Bbb P}^1.
$$ 

2. The complement $\overline {\cal M}_{0,5} - {\cal M}_{0,5}$ is a union of $10$ projective lines. 
A choice of dihedral order of the points $(x_1, ..., x_5)$ provides a splitting of 
them into two pentagons. 
The Stasheff divisor $B_2$ is one of them. See \cite{G} for a 
more elaborate discussion of this example. 

\vskip 3mm

We summarise this discussion in the following Proposition. 
\bp \la{compact2} 

i) The space
 of real positive points ${\cal X}_{A_n}(\R_{>0})$ is isomorphic to ${\cal M}^0_{0,n+3}(\R)$. 

ii) 
The space of the real positive points of $\widehat {\cal X}_{A_n}$ 
is isomorphic to the Stasheff polytope. 

iii)
The cluster modular group of ${\cal X}_{A_n}$
 is identified with $\Z/(n+3)\Z$ if $n>1$ and $\Z/2\Z$ if $n=1$. Its generator acts on 
$\overline {\cal M}_{0,n+3}$ as the cyclic shift
\be \la{n+3}
(x_1,x_2,  ..., x_{n+3}) \lms  (x_2, ..., x_{n+3}, x_1). 
\ee
\ep

\vskip 3mm
Similarly, there is a scheme $\widehat {\cal X}_{\Delta}$ of finite type over $\Z$, the special completion 
of the cluster Poisson variety related to a root system $\Delta$. 
Recall the generalized Stasheff polytope $S_\Delta$ assigned to a root system $\Delta$ \cite{FZII}.  
\bp \la{compactdel}
The generalized Stasheff polytope $S_\Delta$ 
is a combinatorial skeleton of the scheme $\widehat {\cal X}_{\Delta}$. Precisely, 
 $\widehat {\cal X}_{\Delta}$ 
has a stratification whose codimension $k$ strata match  
 the codimension $k$ faces of the polytope $S_\Delta$: 
The latter are products of the generalized Stasheff polytopes $S_{\Delta'}$, and the former 
are products of the corresponding schemes $\widehat {\cal X}_{\Delta'}$. 
One has $\widehat {\cal X}_{\Delta}(\R_{>0}) = S_\Delta$. 
\ep

\subsection{Convexity in the framework of positive spaces.} \la{sec2.4}
\paragraph{Convex subsets of the space of tropical points.} Let ${\cal X}$ be a positive space. 
We define {\it convex subsets} of the space of the tropical points  ${\cal X}({\Bbb A}^t)$, 
where ${\mathbb A}$ is either $\Z$, or $\Q$, or $\R$. 

Let $V$ be a real vector space. A convex domain in $V$ 
is an intersection, possibly infinite, of half spaces given by the inequalities $l(x) \leq c$ where 
$l(x)$ are linear functions on $V$, and $c$ are  constants. 
We look at this as follows. A split torus ${\rm T}$ is the simplest basic example of a positive space. 
The space ${\rm T}(\Z^t)$ of its integral tropical points is a free abelian 
group given by the characters of ${\rm T}$. 
One has  ${\rm T}(\Q^t)= {\rm T}(\Z^t)\otimes \Q$. 
Linear functions with integral coefficients are nothing else but the 
tropicalizations of the positive regular functions on 
the torus ${\rm T}$. We
 generalise this by replacing the torus ${\rm T}$ by an arbitrary positive space ${\cal X}$.

\vskip 3mm
Denote by ${\Bbb L}_+({\cal X})$ the semiring of positive regular functions on 
${\cal X}$. It is the set of rational functions which 
are Laurent polynomials with  positive 
integral coefficients in every coordinate system of the atlas. 
Denote by ${\Bbb E}({\cal X})$ the set of extremal elements of the cone ${\Bbb L}_+({\cal X})$, 
i.e. the elements which can not be presented as a sum of two non-zero 
elements of the cone. 

A function $E \in {\Bbb E}({\cal X})$ and a constant $c_E \in \Q\cup \infty$ give rise to a basic 
convex subset 
$$
C(E, c_E) = \{ x \in {\cal X}(\Q^t) ~|~ E^t(x) \leq c_E\}.
$$
Convex subsets of  ${\cal X}(\Q^t)$ are defined as intersections of 
such subsets. If $c_E = \infty$, 
the corresponding subset is ${\cal X}(\Q^t)$. So if all 
but finitely many $c_E$'s are infinite, 
we get a finite intersection of basic convex subsets. 
We extend the addition in $\Q$ to an addition in $\Q\cup \infty$ by setting $a + \infty=\infty$. 
One defines similarly convex subsets of ${\cal X}(\mathbb A^t)$:

\bd
i) A convex subset $A \subset {\cal X}(\mathbb A^t)$ is  defined by inequalities
$$
A:= \{x \in {\cal X}(\mathbb A^t) ~|~E^t(x) \leq a_E ~
\mbox{for each} ~E \in  {\Bbb E}({\cal X}), \quad a_E \in \mathbb A\cup \infty\}.
$$

ii) Given  $n$ convex subsets 
$$
A_i= \{x \in {\cal X}(\mathbb  A^t) ~|~ E^t(x) \leq a^{(i)}_E ~\mbox{for each} ~E \in {\Bbb E}({\cal X})\},
$$
their Minkowski sum  $A_1 \ast \ldots \ast A_n$ is defined by 
$$
A_1 \ast \ldots \ast A_n := \{x \in {\cal X}(\mathbb A^t) ~|~ E^t(x) \leq a^{(1)}_E+ \ldots + a^{(n)}_E ~\mbox{for each} ~E \in {\Bbb E}({\cal X})\}.
$$
\ed

Each coordinate system identifies the set ${\cal X}(\mathbb A^t)$ with $\mathbb A^{{\rm dim}{\cal X}}$. 
Our convex subsets 
are identified with certain convex (in the usual sense) subsets of $\mathbb A^{{\rm dim}{\cal X}}$.

\paragraph{Spherical convex subsets.} We start again with a motivation. Let $V$ be a rational vector space and ${\mathbb S}(V)$ the sphere of the rays 
in $V$. Convex domains in ${\mathbb S}(V)$ are projectivisation of convex cones 
in $V$. A convex cone in $V$ is an intersection of 
half spaces defined by inequalities $l(v)\leq 0$, where $l$ are linear functions on $V$. 
We conclude that the semiring of positive integral Laurent polynomials 
determines the classical convex structure of the sphere ${\mathbb S}(V)$. 
The real spherical tropical space for a split torus ${\rm T}$ is nothing else but the 
sphere ${\mathbb S}(V)$. This suggests the following generalization.

\bd Let ${\cal X}$ be a positive space. 

i) {\rm Convex subsets}  of the tropical space 
${\cal S}{\cal X}(\mathbb A^t)$  are the subsets $$
S_F = \{ x \in {\cal S}{\cal X}(\Q^t) ~|~ F^t(x) \leq 0\}, \qquad 
F \in {\Bbb L}_+({\cal X}).
$$

ii) The Minkowski sum of convex subsets $S_{F_1} \ast S_{F_2}$ is given by 
$$
S_{F_1} \ast S_{F_2}:= S_{F_1F_2}.
$$
\ed 

Given two  functions $F_1, F_2 \in {\Bbb L}_+({\cal X})$, we have 
\be \la{csubs}
S_{F_1+F_2} = S_{F_1} \cap S_{F_2}. 
\ee
Thus finite intersections of 
the convex subsets in ${\cal X}(\Q^t)$ are convex subsets. This leads to

\bl 
The spherical convex subsets form a semiring with the addition given by the intersection of the 
subsets, and the multiplication given by the Minkowski sum. The map 
$F \lra S_F$ is a morphism of the semiring ${\mathbb L}_+({\cal X})$ 
to the semiring of convex subsets of ${\cal X}$. 
\el

Finally, we can allow infinite intersections of the basic spherical convex subsets.

\end{document}